\documentclass[a4paper,12pt]{amsart}
\usepackage{amsfonts}
\usepackage{amsmath,latexsym,amscd}

\title{Projective 3-folds of general type with $\chi=1$}
\author{Meng Chen and Lei Zhu}
\address{School of Mathematical Sciences, Fudan University, Shanghai, 200433,
People's Republic of China} \email{mchen@fudan.edu.cn}
\email{051018003@fudan.edu.cn}
\thanks{The first author was supported by the Program for New Century Excellent
Talents in University (\#NCET-05-0358) and the National Outstanding
Young Scientist Foundation (\#10625103). The second author was
supported by Graduate Students' Innovation Projects (EYH5928004)}

\newcommand{\bQ}{{\mathbb Q}}

\newcommand{\roundup}[1]{\ulcorner{#1}\urcorner}

\newtheorem{thm}{Theorem}[section]
\newtheorem{lem}[thm]{Lemma}
\newtheorem{cor}[thm]{Corollary}
\newtheorem{prop}[thm]{Proposition}
\newtheorem{claim}[thm]{Claim}
\theoremstyle{definition}
\newtheorem{defn}[thm]{Definition}
\newtheorem{setup}[thm]{}

\newtheorem{rem}[thm]{Remark}
\theoremstyle{remark}

\begin{document}
\begin{abstract} There are many examples of 3-folds of general type
with $\chi(\mathcal {O})=1$ found by Fletcher and Reid about twenty
years ago. Fletcher has ever proved $P_{12}(X)\ge 1$ and
$P_{24}(X)\ge 2$ for all minimal 3-folds $X$ of general type with
$\chi(\mathcal {O}_X)=1$. In this paper, we improve on Fletcher's
method. Our main result is that $\varphi_m$ is birational onto its
image for all $m\geq 63$. To prove this we will show $P_{m}\ge 1$
for all $m\ge 14$ and $P_{2l+18}\geq 3$ for all $l\geq 0$.
\end{abstract}
\maketitle
\pagestyle{myheadings} \markboth{\hfill M. Chen and L.
Zhu\hfill}{\hfill Projective 3-folds of general type\hfill}

\section{\bf Introduction}

To classify algebraic varieties is one of the main goals of
algebraic geometry. In this paper we are concerned with the explicit
algebraic geometry of complex projective 3-folds of general type.

Let $V$ be a smooth projective 3-fold of general type. Let $X$ be a
minimal model of $V$. Denote by $\varphi_m$ the m-th pluricanonical
map. A classic problem is to see when $\varphi_m$ is birational
onto its image. Recently a remarkable theorem by Tsuji
\cite{Tsuji}, Hacon-M$^{\text{\rm c}}$Kernan \cite{H-M} and
Takayama \cite{Tak} says that there is a universal constant $r_3$
such that $\varphi_m$ is birational for all $m\ge r_3$ and for
arbitrary 3-folds of general type. A very new result by J. A. Chen
and the first author in \cite{explicit} shows that one may take
$r_3=77$.

There have been some concrete known bounds on $r_3$ already. For
example, $r_3\le 5$ (sharp) if $X$ is Gorenstein by J. A. Chen, M.
Chen, D.-Q. Zhang \cite{Crelle}; $r_3\le 8$ (sharp) if either
$q(X)>0$ by J. A. Chen, C. D. Hacon \cite{C-H} or $p_g(X)\ge 2$ by
M. Chen \cite{IJM}; $r_3\le 14$ (sharp) if $\chi(\mathcal {O}_X)\le
0$ by M. Chen, K. Zuo \cite{C-Z}. It is natural to study a 3-fold
with $\chi(\mathcal {O})\ge 1$.

First we treat a general 3-fold and prove the following:

\begin{thm}\label{m0} Let $V$ be a nonsingular
projective $3$-fold of general type. Assume $P_{m_1}(V)\ge 2$ and
$P_m(V)\ge 1$ for all $m\ge m_0\ge 2$. Then the pluricanonical map
$\varphi_m$ is birational for all $m\ge \text{max}\{m_0+4m_1+2,
5m_1+4\}$.
\end{thm}

Theorem \ref{m0} has improved Koll\'ar's Corollary 4.8 in \cite{Ko}
and Theorem 0.1 of \cite{JPAA}.

In the second part we prove the following:

\begin{thm}\label{main} Let $V$ be a nonsingular projective 3-fold of
general type with $\chi(\mathcal {O}_V)=1$. Then

(i) $P_m(V):=h^0(V, mK_V)>0$ for all $m\ge 14$;

(ii) $P_{18+2l}(V)\ge 3$ for all integer $l\ge 0$;

(iii) $\varphi_m$ is birational onto its image for all $m\ge 63$.
\end{thm}

Theorem \ref{main} has improved Iano-Fletcher's results in
\cite{Flt}.

Throughout our paper the symbol $\equiv$ stands for the numerical
equivalence of divisors, whereas $\sim$ denotes the linear
equivalence and $=_{\mathbb{Q}}$ denotes the $\mathbb{Q}$-linear
equivalence.

\section{\bf Pluricanonical systems}

In this section we are going to treat a general 3-fold of general
type. By the 3-dimensional MMP (see \cite{K-M, KMM, Reid83} for
instance) we may consider a minimal 3-fold $X$ of general type with
$\mathbb{Q}$-factorial terminal singularities.

\begin{setup}\label{J}{\bf Assumption}. Assume that, on a smooth model
$V_0$ of $X$, there is an effective
divisor $\Gamma\le m_1K_{V_0}$ with $n_{\Gamma}:=h^0(V_0, \mathcal
{O}_{V_0}(\Gamma))\ge 2$. Naturally $P_{m_1}\ge 2$. We would like to
study the rational map $\varphi_{|\Gamma|}$. A very special
situation is $\Gamma= m_1K_{V_0}$, meanwhile $\varphi_{|\Gamma|}$ is
nothing but the $m_1$-canonical map.
\end{setup}

\begin{setup}\label{setup}{\bf Set up.}

First we fix an effective Weil divisor $K_{m_1}\sim m_1K_X$. Take
successive blow-ups $\pi: X'\rightarrow X$ (along nonsingular
centers), which exists by Hironaka's big theorem, such that:

(i) $X'$ is smooth;

(ii) there is a birational morphism $\pi_{\Gamma}:X'\rightarrow
V_0$;

(iii) the movable part $M_{\Gamma}$ of $|\pi^*_{\Gamma}(\Gamma)|$ is
base point free;

(iii) the support of $\pi^*(K_{m_1}) \cup\pi^*_{\Gamma}(\Gamma)$ is
of simple normal crossings.

Denote by $g$ the composition $\varphi_{\Gamma}\circ\pi_{\Gamma}$.
So $g: X'\longrightarrow W'\subseteq{\mathbb P}^{n_{\Gamma}-1}$ is
a morphism. Let $X'\overset{f}\longrightarrow
B\overset{s}\longrightarrow W'$ be the Stein factorization of $g$.
We have the following commutative diagram:\medskip

\begin{picture}(50,80) \put(100,0){$V_0$} \put(100,60){$X'$}
\put(170,0){$W'$} \put(170,60){$B$} \put(112,65){\vector(1,0){53}}
\put(106,55){\vector(0,-1){41}} \put(175,55){\vector(0,-1){43}}
\put(114,58){\vector(1,-1){49}} \multiput(112,2.6)(5,0){11}{-}
\put(162,5){\vector(1,0){4}} \put(133,70){$f$} \put(180,30){$s$}
\put(92,30){$\pi_{\Gamma}$}
\put(135,-8){$\varphi_{\Gamma}$}\put(136,40){$g$}
\end{picture}
\bigskip

Denote by $M_{k}$ the movable part of $|kK_{X'}|$ for any positive
integer $k>0$. We may write $m_1K_{X'}=_{\mathbb
Q}\pi^*(m_1K_X)+E_{\pi, m_1}=M_{m_1}+Z_{m_1},$ where $M_{m_1}$ is
the movable part of $|m_1K_{X'}|$, $Z_{m_1}$ the fixed part and
$E_{\pi, m_1}$ an effective ${\mathbb Q}$-divisor which is a
${\mathbb Q}$-sum of distinct exceptional divisors. By $K_{X'}-
\frac{1}{m_1}E_{\pi, m_1}$, we mean $\pi^*(K_X)$. So, whenever we
take the round up of $m\pi^*(K_X)$, we always have
$\roundup{m\pi^*(K_X)}\le mK_{X'}$ for all positive numbers $m$.
Since $M_{\Gamma}\le M_{m_1}\le \pi^*(m_1K_X)$, we can write
$\pi^*(m_1K_X)=M_{\Gamma}+E_{\Gamma}'$ where $E_{\Gamma}'$ is an
effective $\mathbb{Q}$-divisor. Set $d:=\dim (B)$. Denote by $S$ a
generic irreducible element (see Definition \ref{generic}) of
$|M_{\Gamma}|$. Then $S$ is a smooth projective surface of general
type. When $d=1$, one can write $M_{\Gamma}\equiv a_{\Gamma} S$
where $a_{\Gamma}\ge n_{\Gamma}-1$.
\end{setup}

\begin{defn}\label{generic}
 Assume that a complete linear system $|M'|$ is movable on an
 arbitrary variety $V$. {\it
A generic irreducible element} $S'$ of $|M'|$ is defined to be a
generic irreducible component in a general member of $|M'|$.
Clearly $S'\sim M'$ only when $|M'|$ is not composed with a pencil.
\end{defn}

Before proving the main result we build a technical, but a quite
useful theorem which is a generalized form of Theorem 2.6 in
\cite{IJM}.


\begin{thm}\label{Key} Let $X$ be a minimal
projective $3$-fold of general type with $\mathbb{Q}$-factorial
terminal singularities. Assume that, on a smooth model $V_0$ of $X$,
there is an effective divisor $\Gamma\le m_1K_{V_0}$ with
$n_{\Gamma}:=h^0(V_0, \mathcal {O}_{V_0}(\Gamma))\ge 2$. Keep the
same notation as in $\ref{setup}$ above. One has a fibration $f:
X'\longrightarrow B$ induced by $\varphi_{\Gamma}$. Denote by $S$ a
generic irreducible element of $|M_{\Gamma}|$. Assume that, on the
smooth surface $S$, there is a movable linear system $|G|$ (not
necessarily base point free). Let $C$ be a generic irreducible
element of $|G|$ (so $C$ can be singular). Set
$\xi:=(\pi^*(K_X)\cdot C)_{X'}\in \mathbb{Q}$ and
$$p:=\begin{cases} 1  &\text{if}\ \dim (B)\ge 2\\
a_{\Gamma} (\text{see}\ \ref{setup}\ \text{for the definition})
&\text{otherwise.}
\end{cases}$$
Then the inequality $m\xi\ge 2g(C)-2+\alpha_0$ (where $g(C)$ is the
geometric genus of $C$) holds under the assumptions $(1)$ and $(2)$
below. Furthermore $\varphi_m$ of $X$ is birational onto its image
under the assumptions $(1), (2)', (3)$ and $(4)$ below.

Assumptions, for a positive integer $m$:
\begin{itemize}

\item[(1)\ ]
There is a rational number $\beta>0$ such that
${\pi}^*(K_X)|_{S}-\beta C$ is numerically equivalent to an
effective ${\mathbb Q}$-divisor; and set
$\alpha:=(m-1-\frac{m_1}{p}-\frac{1}{\beta})\xi$ and
$\alpha_0:=\roundup{\alpha}$.
\item[(2)\ ]
The inequality $\alpha > 1$ holds.

\item[(2)']
Either $\alpha>2$ or $\alpha>1$ and $C$ is non-hyper-elliptic.

\item[(3)\ ] The linear system $|mK_{X'}|$ separates different
generic irreducible elements of $|M_{\Gamma}|$ (namely,
$\Phi_{|mK_{X'}|}(S')\ne \Phi_{|mK_{X'}|}(S'')$ for two different
generic irreducible elements $S'$, $S''$ of $|M_{\Gamma}|$).

\item[(4)\ ] The linear system $|mK_{X'}||_S$ on $S$ (as a sub-linear system
of $|mK_{X'}|_S|$) separates different generic irreducible elements
of $|G|$. Or sufficiently, the complete linear system
$$|K_{S}
+\roundup{(m-1)\pi^*(K_X)-S-\frac{1}{p}E_{\Gamma}'}|_{S}|$$
separates different generic irreducible elements of $|G|$.
\end{itemize}

\end{thm}

\begin{proof} We first prove the birationality of $\varphi_m$.

Condition (3) says that the linear system $|mK_{X'}|$ separates
different irreducible elements of $|M_{\Gamma}|$. By the
birationality principle (P1) and (P2) of \cite{MPCPS}, it is
sufficient to prove that the linear system $|mK_{X'}||_S$ gives a
birational map on a generic irreducible element $S$ of
$|M_{\Gamma}|$. Condition (4) says that $|mK_{X'}||_S$ on $S$
separates different generic irreducible elements of $|G|$. Again by
the birationality principle it suffices to prove the birationality
of $\Phi_{|mK_{X'}|}|_C$ where $C$ is a generic irreducible element
of $|G|$. In fact, we consider a smaller linear system than
$|mK_{X'}|$. we consider the sub-system
$$|K_{X'}+\roundup{(m-1)\pi^*(K_X)-
\frac{1}{p}E_{\Gamma}'}|\subset |mK_{X'}|.$$ Noting that
$(m-1)\pi^*(K_X)-\frac{1}{p}E_{\Gamma}'-S\equiv
(m-1-\frac{m_1}{p})\pi^*(K_X)$ is nef and big under the assumptions
(1) or (2), the vanishing theorem gives a surjective map
$$H^0(X',K_{X'}+\roundup{(m-1)\pi^*(K_X)-\frac{1}{p}E_{\Gamma}'})$$
$$\longrightarrow  H^0(S,
K_{S}+\roundup{(m-1)\pi^*(K_X)-S-\frac{1}{p}E_{\Gamma}'}|_{S}).\eqno{(2.1)}
$$
Note that
$|K_{X'}+\roundup{(m-1)\pi^*(K_X)-\frac{1}{p}E_{\Gamma}'}|\subset
|mK_{X'}|$. It suffices to study
$|K_{S}+\roundup{(m-1)\pi^*(K_X)-S-\frac{1}{p}E_{\Gamma}'}|_{S}|.$
We go on studying the restriction to $C$. By assumption (1), there
is an effective ${\mathbb Q}$-divisor $H$ on $S$ such that
$$\frac{1}{\beta}\pi^*(K_X)|_{S}\equiv C+H.$$
If $|G|$ is not base point free, we can take a birational
modification $\nu:\overline{S}\longrightarrow S$ (otherwise set
$\nu$ to be an identity) such that $\overline{S}$ is smooth and that
the movable part $|\overline{G}|$ of $|\nu^*(G)|$ is base point
free. Denote by $\overline{C}$ a generic irreducible element of
$|\overline{G}|$. Then $\overline{C}$ is smooth. We study the linear
system
$$|K_{\overline{S}}+\roundup{\nu^*(((m-1)\pi^*(K_X)-S-
\frac{1}{p}E_{\Gamma}')|_{S})}|$$ which is smaller than the linear
system
$$|K_{\overline{S}}+\nu^*(\roundup{((m-1)\pi^*(K_X)-S-
\frac{1}{p}E_{\Gamma}')|_{S}})|.$$ The last linear system has the
same birationality as that of
$|K_{S}+\roundup{((m-1)\pi^*(K_X)-S-\frac{1}{p}E_{\Gamma}')|_{S}}|$
on the surface $S$. We write $\nu^*(C)=\overline{C}+\overline{E}$
where $\overline{E}$ is exceptional and effective. Whenever
$m-1-\frac{m_1}{p}-\frac{1}{\beta}>0$, the ${\mathbb Q}$-divisor
$$\nu^*(((m-1)\pi^*(K_X)-S-\frac{1}{p}E_{\Gamma}')|_{S}-H)-\overline{C}-\overline{E}$$ $$\equiv
(m-1-\frac{m_1}{p}-\frac{1}{\beta})\pi^*(K_X)|_S$$ is nef and big.
The Kawamata-Viehweg vanishing gives the surjective map
$$H^0(\overline{S},
K_{\overline{S}}+\roundup{\nu^*(((m-1)\pi^*(K_X)-S-\frac{1}{p}E_{\Gamma}')|_{S}-H)-\overline{E}})$$
$$\longrightarrow  H^0(\overline{C}, K_{\overline{C}} + D)\hskip6.8cm \eqno (2.2)$$
where
$D:=\roundup{\nu^*(((m-1)\pi^*(K_X)-S-\frac{1}{p}E_{\Gamma}')|_{S}-C-H)-\overline{E}}|_{\overline{C}}$
is a divisor on $\overline{C}$. Noting that $\overline{C}$ is nef on
$\overline{S}$, we have
$$\deg(D)\geq \nu^*(((m-1)\pi^*(K_X)-S-\frac{1}{p}E_{\Gamma}')|_{S}-C-H)\cdot
\overline{C}$$
$$=(((m-1)\pi^*(K_X)-S-\frac{1}{p}E_{\Gamma}')|_{S}-C-H)\cdot C$$ $$=(m-1-\frac{m_{1}}{p}-\frac{1}{\beta})\pi^*(K_{X})|_{S}\cdot C=\alpha$$ and thus $\deg(D)\geq \alpha_0$. Whenever $\deg(D)\ge 3$
or $\deg(D)\ge 2$ and $\overline{C}$ is non-hyper-elliptic,
$$|K_{\overline{S}}+\roundup{\nu^*(((m-1)\pi^*(K_X)-S-\frac{1}{p}E_{\Gamma}')|_{S}-H)-\overline{E}}||_{\overline{C}}$$
gives a birational map. Thus
$|K_{S}+\roundup{(m-1)\pi^*(K_X)-S-\frac{1}{p}E_{\Gamma}'}|_{S}|$
gives a birational map and so $\varphi_m$ of $X$ is birational.

Finally we show the inequality for $\xi$. Whenever we have
$\deg(D)\ge \alpha\ge 2$, $|K_{\overline{C}}+D|$ is base point free
by the curve theory. Denote by $|N_m|$ the movable part of
$|K_{S}+\roundup{((m-1)\pi^*(K_X)-S-\frac{1}{p}E_{m_1}')|_{S}-H}|$.
Applying Lemma 2.7 of \cite{MPCPS} to surjective maps (B1) and (B2),
we have
$$m\pi^*(K_X)|_{S}\ge N_m\ \ \text{and}\ \ (N_m\cdot C)_{S}=\nu^*(N_m)\cdot \overline{C}\ge 2g(C)-2
+\deg (D).$$ Note that the above inequality holds without conditions
(3) or (4). We are done.
\end{proof}


\begin{thm}\label{d3} Let $X$ be a minimal
projective $3$-fold of general type with $\mathbb{Q}$-factorial
terminal singularities. Assume that, on a smooth model $V_0$ of $X$,
there is an effective divisor $\Gamma\le m_1K_{V_0}$ with
$n_{\Gamma}:=h^0(V_0, \mathcal {O}_{V_0}(\Gamma))\ge 2$. Assume
$P_{m}\ge 1$ for all $m\ge m_0\ge 2$. Keep the same notation as in
\ref{setup}. If $d=3$, Then $\varphi_m$ is birational for all $m\ge
\text{max}\{m_0+m_1, 3m_1+2\}$.
\end{thm}
\begin{proof}Recall that we have $p=1$.

Take an integer $m\ge m_0+m_1$. Since $mK_{X'}\ge M_{\Gamma}$ and
that $|M_{\Gamma}|$ is not composed  with a pencil, $|mK_{X'}|$ can
separate different $S$. Theorem \ref{Key}(3) is satisfied. On the
surface $S$, we take $G:=S|_S$. Then $|G|$ is not composed of a
pencil. Since we have
$$K_{S}+\roundup{(m-1)\pi^*(K_X)-S-E_{\Gamma}'}|_{S}\ge
S|_S,$$ the exact sequence (B1) shows that Theorem \ref{Key}(4) is
satisfied.

Because $m_1\pi^*(K_X)|_S\ge C$ where $C\in |G|$ is a general
member, we can take $\beta=\frac{1}{m_1}$. So Theorem \ref{Key}(1)
is satisfied. On a generic irreducible element $S$ of
$|M_{\Gamma}|$, we have a linear system $|G|$ which is not composed
of a pencil and is base point free. So $S^3=S|_S\cdot S|_S=C^2\ge
2$. Thus
$$2g(C)-2=(K_S+C)\cdot C\ge (\pi^*(K_X)|_S+2S|_S)\cdot C>4,$$
which means $g(C)\ge 4$.

If we take a sufficiently big $m$, then
$(m-1-\frac{m_1}{p}-\frac{1}{\beta})\xi$ will be big enough.
Theorem \ref{Key} gives $\xi\ge \frac{6}{2m_1+1}$. Take $m\ge
3m_1+2$. Then $\alpha=(m-1-\frac{m_1}{p}-\frac{1}{\beta})\xi\ge
\frac{6m_1+6}{2m_1+1}>2$. (Theorem \ref{Key} gives $\xi\ge
\frac{9}{3m_1+2}$.)

Therefore one sees by Theorem \ref{Key} that $\varphi_m$ is
birational for all $m\ge \text{max}\{m_0+m_1, 3m_1+2\}$.

One can even get better bound of m whenever $m_1$ is big. For
instance, if $m_1\ge 11$, take $m\ge 3m_1-2$. Then
$\alpha=(m-2m_1-1)\xi\ge \frac{9m_1-27}{3m_1+2}>2$. So $\varphi_m$
is birational for all $m\ge \text{max}\{m_0+m_1, 3m_1-2\}$ and
$m_1\geq 11$.
\end{proof}


\begin{thm}\label{d2} Let $X$ be a minimal
projective $3$-fold of general type with $\mathbb{Q}$-factorial
terminal singularities. Assume that, on a smooth model $V_0$ of $X$,
there is an effective divisor $\Gamma\le m_1K_{V_0}$ with
$n_{\Gamma}:=h^0(V_0, \mathcal {O}_{V_0}(\Gamma))\ge 2$. Assume
$P_m\ge 1$ for all $m\ge m_0\ge 2$. Keep the same notation as in
\ref{setup}. If $d=2$, then $\varphi_m$ is birational for all $m\ge
\text{max}\{m_0+2m_1, 4m_1+2\}$. Furthermore $\varphi_m$ is
birational for all $m\ge \text{max}\{m_0+2m_1, 4m_1-9\}$ and for
$m_1\geq 18$.
\end{thm}

\begin{proof} Recall that we have $p=1$. Take an integer $m\ge m_0+2m_1$. Since $mK_{X'}\ge M_{\Gamma}$ and
that $|M_{\Gamma}|$ is not composed with a pencil, $|mK_{X'}|$ can
separate different $S$. Theorem \ref{Key}(3) is satisfied. On the
surface $S$, we take $G:=S|_S$. Different from the case $d=3$, $|G|$
is composed with a pencil of curves. If $|G|$ is composed of a
rational pencil, then since we have
$$K_{S}+\roundup{(m-1)\pi^*(K_X)-S-E_{\Gamma}'}|_{S}\ge
S|_S,$$ the exact sequence (2.1) shows that $|mK_{X'}||_S$ can
separate different generic irreducible elements of $|G|$. Otherwise
$G\equiv tC$ where $t>1$ and $C$ is a generic irreducible element
of $|G|$. Noting that $m_1\pi^*(K_X)|_S\equiv
tC+{E_{\Gamma}'}|_{S}$, one sees that
\begin{eqnarray*}
&&(m-1)\pi^*(K_X)|_S-S|_S-(1+\frac{2}{t})E_{\Gamma}'|_S-C_1-C_2\\
&\equiv & (m-1-m_1-\frac{2m_1}{t})\pi^*(K_X)|_S
\end{eqnarray*}
is
nef and big, where $C_1$ and $C_2$ are different generic irreducible
elements of $|G|$. Thus the Kawamata-Viehweg vanishing theorem
(\cite{Kav,V}) gives the surjective map:
\begin{eqnarray*}
&& H^0(S,K_{S}+\roundup{((m-1)\pi^*(K_X)-S-
(1+\frac{2m_1}{t})E_{\Gamma}')|_{S}})\\
&\longrightarrow & H^0(C_1, K_{C_1}+D_1)\oplus H^0(C_2,K_{C_2}+
D_2)\longrightarrow 0
\end{eqnarray*}
where $\deg(D_i)>0$ for $i=1,2$. Because $C_i$ is a curve of genus
$\ge 2$, $H^0(C_i,K_{C_i}+D_i)\neq 0$ for $i=1,2$. So $|mK_{X'}||_S$
can separate different generic irreducible elements of $|G|$. In a
word, Theorem \ref{Key}(4) is satisfied for all $m\ge m_0+2m_1$.

Because $m_1\pi^*(K_X)|_S\ge C$ where $C$ is a generic irreducible
of $|G|$, we can take $\beta=\frac{1}{m_1}$. So Theorem \ref{Key}(1)
is satisfied. On a generic irreducible element $S$ of
$|M_{\Gamma}|$, we have a linear system $|G|$ which is composed of a
pencil of curves. Because $g(C)\ge 2$, one has $K_S\cdot C+C^2\ge
2$.

If we take a sufficiently big $m$, then
$(m-1-\frac{m_1}{p}-\frac{1}{\beta})\xi$ will be big enough. Theorem
\ref{Key} gives $\xi\ge \frac{2}{2m_1+1}$. Take $m=4m_1+3$. Then
$\alpha=(m-\frac{m_1}{p}-\frac{1}{\beta})\xi\ge
\frac{4m_1+4}{2m_1+1}>2$. Theorem \ref{Key} gives $\xi\ge
\frac{5}{4m_1+3}$. Take $m=4m_1+2$. Then $\alpha\ge
\frac{10m_1+5}{4m_1+3}>2$.

Therefore Theorem \ref{Key} says that $\varphi_m$ is birational
whenever $m\ge \text{max}\{m_0+2m_1, 4m_1+2\}$. Also one gets
$\xi\geq \frac{5}{4m_1+2}$.

Now assume $m_1\ge 10$. Take $m\geq 4m_1-2$ and then $\alpha>2$. So
$\varphi_m$ is birational for all $m\ge \text{max}\{m_0+2m_1,
4m_1-2\}$; Assume $m_1\geq 16$ and take $m\geq 4m_1-6$. Then
$\alpha>2$ and $\xi\geq \frac{5}{4m_1-6}$. Take $m\geq 4m_1-8$. Then
$\alpha>2$ and $\xi\geq \frac{5}{4m_1-8}$; Assume $m_1\geq 18$ and
take $m\geq 4m_1-9$. Then $\alpha>2$. In a word, we have seen that
$\varphi_m$ is birational whenever $m\ge \text{max}\{m_0+2m_1,
4m_1-9\}$ provided $m_1\geq 18$.
\end{proof}

Now we begin to study the case $d=1$. Though similar lemmas has
already been established in several papers of the first author, we
include a more general one here for the convenience to future
applications.

\begin{lem}\label{b>0} Let $X$ be an arbitrary minimal 3-fold with
${\mathbb Q}$-factorial terminal singularities. Let
$\pi:X'\longrightarrow X$ be a smooth birational modification.
Assume that $f:X'\longrightarrow B$ is an arbitrary fibration onto a
smooth curve $B$ with $g(B)>0$. Denote by $F$ a general fiber of
$f$. Then $\pi^*(K_X)|_F\sim\sigma^*(K_{F_0})$ where
$\sigma:F\longrightarrow F_0$ is the blow down onto the smooth
minimal model.
\end{lem}
\begin{proof} We shall use the idea of Lemma 14 in Kawamata's paper
\cite{KA}. By Shokurov's theorem in \cite{Sho}, each fiber of
$\pi:X'\longrightarrow X$ is rationally chain connected. Therefore,
$f(\pi^{-1}(x))$ is a point for all $x\in X$. Considering the image
$G\subset (X \times B)$ of $X'$ via the morphism $(\pi\times f)\circ
\triangle_{X'}$ where $\triangle_{X'}$ is the diagonal map
$X'\longrightarrow X'\times X'$, one knows that $G$ is a projective
variety. Let $g_1:G\longrightarrow X$ and $g_2:G\longrightarrow B$
be two projections. Since $g_1$ is a projective morphism and even a
bijective map, $g_1$ must be both a finite morphism of degree 1 and
a birational morphism. Since $X$ is normal, $g_1$ must be an
isomorphism. So $f$ factors as $f_1 \circ \pi$ where $f_1:=g_2\circ
g_1^{-1} : X \rightarrow B$ is a well defined morphisms. In
particular, a general fiber $F_0$ of $f_1$ must be smooth minimal.
So it is clear that $\pi^*(K_X)|_F\sim\sigma^*(K_{F_0})$ where
$\sigma$ is nothing but $\pi|_F$.
\end{proof}

The following lemma shows a way to find a suitable $\beta$ in
Theorem \ref{Key}. We admit that it has already appeared as a weaker
form in several papers of the first author.

\begin{lem}\label{beta} Let $X$ be a minimal
projective $3$-fold of general type with $\mathbb{Q}$-factorial
terminal singularities. Assume that, on a smooth model $V_0$ of $X$,
there is an effective divisor $\Gamma\le m_1K_{V_0}$ with
$n_{\Gamma}:=h^0(V_0, \mathcal {O}_{V_0}(\Gamma))\ge 2$. Assume
$P_m\ge 1$ for all $m\ge m_0\ge 2$. Keep the same notation as in
\ref{setup}. Assume $B=\mathbb{P}^1$. Let $f:X'\longrightarrow
\mathbb{P}^1$ be an induced fibration of $\varphi_{|\Gamma|}$.
Denote by $F:=S$ a general fiber of $f$. Then one can find a
sequence of rational numbers $\{\beta_n\}$ with $\lim_{n\mapsto
+\infty} \beta_n = \frac{p}{m_1+p}$ such that
$\pi^*(K_X)|_F-\beta_n\sigma^*(K_{F_0})$ is ${\bQ}$-linearly
equivalent to an effective ${\bQ}$-divisor $H_n$, where
$\sigma:F\longrightarrow F_0$ is the blow down onto the smooth
minimal model.
\end{lem}

\begin{proof}

One has $\mathcal {O}_{B}(p)\hookrightarrow {f}_*\mathcal
{O}_{X'}(M_{\Gamma})\hookrightarrow {f}_*\omega_{X'}^{m_1}$ and
therefore ${f}_*\omega_{X'/B}^{t_0p}\hookrightarrow
{f}_*\omega_{X'}^{t_0p+2t_0m_1}$ for any big integer $t_0$.

For any positive integer $k$, denote by $M_k$ the movable part of
$|kK_{X'}|$. Note that ${f}_*\omega_{X'/B}^{t_0p}$ is generated by
global sections since it is semi-positive according to E. Viehweg
(\cite{VV}). So any local section can be extended to a global one.
On the other hand, $|t_0p\sigma^*(K_{F_0})|$ is base point free and
is exactly the movable part of $|t_0pK_F|$ by Bombieri \cite{Bom}.
Set $a_0:=t_0p+2t_0m_1$ and $b_0:=t_0p$. Clearly one has the
following relation:
$$a_0\pi^*(K_X)|_F\ge M_{t_0p+2t_0m_1}|_F\ge b_0\sigma^*(K_{F_0}).$$
 This means that there is an
effective $\mathbb{Q}$-divisor $E_0'$ on $F$ such that
$$a_0\pi^*(K_X)|_F=_{\bQ} b_0\sigma^*(K_{F_0})+E_0'.$$
Thus $\pi^*(K_X)|_F =_{\bQ} \frac{p}{p+2m_1}\sigma^*(K_{F_0})+E_0$
with $E_0=\frac{1}{a_0}E_0'$.

Let us consider the case $p\ge 2$.

Assume that we have defined $a_n$ and $b_n$ such that the following
is satisfied with $l = n:$
$$a_{l}\pi^*(K_X)|_F \ge b_{l}\sigma^*(K_{F_0}).$$
We will define $a_{n+1}$ and $b_{n+1}$ inductively such that the
above inequality is satisfied with $l = n+1$. One may assume from
the beginning that $a_n\pi^*(K_X)$ is supported on a divisor with
normal crossings. Then the Kawamata-Viehweg vanishing theorem
implies the surjective map
$$H^0(K_{X'}+\roundup{a_n\pi^*(K_X)}+F)\longrightarrow H^0(F, K_F+
\roundup{a_n\pi^*(K_X)}|_F).$$ One has the relation
\begin{eqnarray*}
|K_{X'}+\roundup{a_n\pi^*(K_X)}+F||_F&=&|K_F+\roundup{a_n\pi^*(K_X)}|_F|\\
&\supset& |K_F+b_n\sigma^*(K_{F_0})|\\
&\supset& |(b_n+1)\sigma^*(K_{F_0})|.
\end{eqnarray*}
Denote by $M_{a_n+1}'$ the movable part of $|(a_n+1)K_{X'}+F|$.
Applying Lemma 2.7 of \cite{MPCPS}, one has $M_{a_n+1}'|_F\ge
(b_n+1)\sigma^*(K_{F_0}).$ We modify our original $\pi$ such that
$|M_{a_n+1}'|$ is base point free. In particular, $M_{a_n+1}'$ is
nef. Since $X$ is of general type $|mK_X|$ gives a birational map
whenever $m$ is big enough. Thus we see that $M_{a_n+1}'$ is big if
we fix a very big $t_0$ in advance.

Now the Kawamata-Viehweg vanishing theorem again gives
\begin{eqnarray*}
|K_{X'}+M_{a_n+1}'+F||_F&=&|K_F+M_{a_n+1}'|_F|\\
&\supset& |K_F+(b_n+1)\sigma^*(K_{F_0})|\\
&\supset& |(b_n+2)\sigma^*(K_{F_0})|.
\end{eqnarray*}

We may repeat the above procedure inductively. Denote by
$M_{a_n+t}'$ the movable part of $|K_{X'}+M_{a_n+t-1}'+F|$ for $t\ge
2$. For the same reason, we may assume $|M_{a_n+t}'|$ to be base
point free. Inductively one has:
$$M_{a_n+t}'|_F\ge (b_n+t)\sigma^*(K_{F_0}).$$
Applying the vanishing theorem once more, we have
\begin{eqnarray*}
|K_{X'}+M_{a_n+t}'+F||_F&=&|K_F+M_{a_n+t}'|_F|\\
&\supset& |K_F+(b_n+t)\sigma^*(K_{F_0})|\\
&\supset& |(b_n+t+1)\sigma^*(K_{F_0})|.
\end{eqnarray*}

Take $t=p-1$. Noting that
$$|K_{X'}+M_{a_n+p-1}'+F|\subset |(a_n+p+m_1)K_{X'}|$$
and applying Lemma 2.7 of \cite{MPCPS} again, one has
$$a_{n+1}\pi^*(K_X)|_F\ge M_{a_n+p+m_1}|_F\ge M'_{a_n+p}|_F\ge
b_{n+1} \sigma^*(K_{F_0}).$$ Here we set $a_{n+1}:=a_n+p+m_1$ and
$b_{n+1}=b_n+p$. Set $\beta_n = \frac{b_{n}}{a_{n}}.$ Clearly
$\lim_{n\mapsto +\infty} \beta_n = \frac{p}{m_1+p}$.

The case $p=1$ can be proved similarly, but with a simpler
induction. We omit the details.
\end{proof}


\begin{thm}\label{d1} Let $X$ be a minimal projective $3$-fold of
general type with $\mathbb{Q}$-factorial terminal singularities.
Assume that, on a smooth model $V_0$ of $X$, there is an effective
divisor $\Gamma\le m_1K_{V_0}$ with $n_{\Gamma}:=\newline h^0(V_0,
\mathcal {O}_{V_0}(\Gamma))\ge 2$. Assume $P_m\ge 1$ for all $m\ge
m_0\ge 2$. Keep the same notation as in \ref{setup}. If $d=1$, then:
\begin{itemize}
\item[(1)]
$\varphi_m$ is birational for all $m\ge \text{max}\{m_0+4m_1+2,
5m_1+4\}$;

\item[(2)] $\varphi_m$ is birational for all $m\ge
\text{max}\{m_0+4m_1+2, 5m_1-2\}$ and for $m_1\ge 14$;

\item[(3)] Whenever $n_{\Gamma}\geq 3$, $\varphi_{m}$ is birational for all $m\ge
\text{max}\{m_0+2m_1+2, 3m_1+4\}$.
\end{itemize}
\end{thm}
\begin{proof} One has an induced fibration $f:X'\longrightarrow B$.
Denote by $F$ a general fiber of $f$.
\medskip

{\bf Case 1}. $b=g(B)>0$.  Under this situation we have
$p=a_{\Gamma}\ge 2$. By Lemma \ref{b>0}, we have $\pi^*(K_X)|_F\sim
\sigma^*(K_{F_0})$. Let $F'$ and $F''$ be two different smooth
fibers of $f$. For all $m\ge m_0+m_1+3$ one has
$$(m-1)\pi^*(K_X)-\frac{2}{p}E_{\Gamma}'-F'-F''\equiv
(m-1-\frac{2m_1}{p})\pi^*(K_X)$$ is nef and big, the
Kawamata-Viehweg vanishing theorem (\cite{Kav, V}) gives a
surjective map:
\begin{eqnarray*}
&& H^0(X',K_{X'}+\roundup{(m-1)\pi^*(K_X)-\frac{2}{p}E_{\Gamma}'})\\
&\longrightarrow& H^0(F',
K_{F'}+\roundup{(m-1)\pi^*(K_X)-\frac{2}{p}E_{\Gamma}'}|_{F'})\oplus\\
&&H^0(F'',
K_{F''}+\roundup{(m-1)\pi^*(K_X)-\frac{2}{p}E_{\Gamma}'}|_{F''}).
\end{eqnarray*}
The last two groups are non-zero because, for instance,
$$K_{F''}+\roundup{(m-1)\pi^*(K_X)-\frac{2}{p}E_{\Gamma}'}|_{F''}\ge
(m_0+3)\sigma''^*(K_{F_0''})>0$$ where $\sigma'':F''\longrightarrow
F_0''$ is the blow down onto the minimal model. Therefore
$|K_{X'}+\roundup{(m-1)\pi^*(K_X)-\frac{2}{p}E_{\Gamma}'}|$ can
separate $F'$ and $F''$ and so can $|mK_{X'}|$. The vanishing
theorem gives another surjective map:
\begin{eqnarray*}
&& H^0(X',K_{X'}+\roundup{(m-1)\pi^*(K_X)-\frac{1}{p}E_{\Gamma}'})\\
&\longrightarrow& H^0(F,
K_{F}+\roundup{(m-1)\pi^*(K_X)-\frac{1}{p}E_{\Gamma}'}|_{F}).
\end{eqnarray*}
Because
$$K_{F}+\roundup{(m-1)\pi^*(K_X)-\frac{1}{p}E_{\Gamma}'}|_{F}\ge
(m_0+3)\sigma^*(K_{F_0})\ge 5\sigma^*(K_{F_0}),$$ one sees that
$|K_{F}+\roundup{(m-1)\pi^*(K_X)-\frac{1}{p}E_{\Gamma}'}|_{F}|$
gives a birational map because $\varphi_{|5K_F|}$ is birational
according to Bombieri. Therefore $\varphi_m$ is birational onto its
image for all $m\ge m_0+m_1+3$.
\medskip

{\bf Case 2}. $b=g(B)=0$. Suppose $m\ge m_0+4m_1+2$. Consider the
linear system
$|K_{X'}+\roundup{(m-1)\pi^*(K_X)-\frac{1}{p}E_{\Gamma}'}|\subset
|mK_{X'}|$. Because
$K_{X'}+\roundup{(m-1)\pi^*(K_X)-\frac{1}{p}E_{\Gamma}'}\ge F$ and
$|M_{\Gamma}|$ is composed with a rational pencil, $|mK_{X'}|$ can
separate different generic irreducible elements of $|M_{\Gamma}|$.
Theorem \ref{Key} (3) is satisfied. Note that $P_2(F)\ge 2$ by
surface theory. On a general fiber $F$, take $G$ to be the movable
part of $|2\sigma^*(K_{F_0})|$. Let $C$ be a generic irreducible
element of $|G|$. According to Theorem 1 of Xiao (\cite{X}), $|G|$
is composed with a pencil of curves if and only if $K_{F_0}^2=1$ and
$p_g(F_0)=0$. Clearly $|G|$ is composed of a rational pencil when
$K_{F_0}^2=1$ and $p_g(F_0)=0$ since $q(F_0)=0$. So it suffices to
show $|mK_{X'}||_F\supset |G|$ in order to verify Theorem
\ref{Key}(4). In fact, one has $\mathcal {O}_{B}(1)\hookrightarrow
f_*\omega_{X'}^{m_1}$. Thus there is the inclusion:
$$f_*\omega_{X'/B}^2\hookrightarrow f_*\omega_{X'}^{4m_1+2}.$$
Since $f_*\omega_{X'/B}^2$ is semi-positive by Viehweg \cite{VV} and
thus generated by global sections, one has
$$m\pi^*(K_X)|_F\ge (4m_1+2)\pi^*(K_X)|_F\ge G.$$
In a word, Theorem \ref{Key} (3) and (4) are satisfied for all $m\ge
m_0+4m_1+2$.

Take a $\beta$, nearby $\frac{1}{2m_1+2}$, by virtue of Lemma
\ref{beta}.

If $|G|$ is composed with a pencil of curves and $g(C)=2$, then
 $\sigma^*(K_{F_0})\cdot C\ge 2$ by Lemma
\ref{10}. This gives $\xi\ge \frac{1}{m_1+1}\sigma^*(K_{F_0})\cdot
C\ge \frac{2}{m_1+1}$ by Lemma 2.8. Take $m\ge 5m_1+3$. Then
$\alpha\ge 2+\frac{2m_1-2}{m_1+1}>2$. Thus Theorem \ref{Key} says
that $\varphi_m$ is birational for all $m\ge \text{max}\{m_0+4m_1+4,
5m_1+3\}$.

Otherwise $|G|$ is composed with a pencil of curves and $g(C)\ge 3$
or $|G|$ is not composed of a pencil of curves. In the later case,
after a necessary birational modification to get the base point
freeness of $|G|$, one sees that $2g(C)-2=K_S\cdot C+C^2\ge 4$.
Again $g(C)\ge 3$. If we take a very large $m$ such that $\alpha$ is
big enough, then Theorem \ref{Key} gives
$$m\xi\ge 2g(C)-2+\roundup{(m-1-\frac{m_1}{p}-\frac{1}{\beta})\xi},$$
which means $\xi\ge \frac{4}{3m_1+3}$. Take $m\ge 5m_1+4$. Then
$\alpha\ge 2+\frac{2m_1-2}{3m_1+3}>2$ whenever $m_1>1$. Therefore
Theorem \ref{Key} says that $\varphi_{m}$ is birational for all
$m\ge \text{max}\{m_0+4m_1+4, 5m_1+4\}$ and for $m_1>1$. The same
statement for the situation $m_1=1$ was proved in \cite{IJM}.

When $m_1$ is big, one can get better bound of $m$. For example,
when $m_1\ge 14$, one sees that $\varphi_m$ is birational for all
$m\ge \text{max}\{m_0+4m_1+2, 5m_1-2\}.$
\medskip

{\bf Case 3}. $n_{\Gamma}\geq 3$. This is more or less  parallel to
Case 2. But since $n_{\Gamma}$ is bigger, we hope to deduce a better
result on the birationality of $\varphi_m$. The case with $b>0$
follows from Case 1. So we may still assume $b=g(B)=0$.

Suppose $m\ge m_0+2m_1+2$. Consider the linear system
$|K_{X'}+\roundup{(m-1)\pi^*(K_X)-\frac{1}{p}E_{\Gamma}'}|\subset
|mK_{X'}|$. Because
$$K_{X'}+\roundup{(m-1)\pi^*(K_X)-\frac{1}{p}E_{\Gamma}'}\geq F$$ and
$|M_{\Gamma}|$ is composed with a rational pencil, $|mK_{X'}|$ can
separate different generic irreducible elements of $|M_{\Gamma}|$.
Theorem \ref{Key} (3) is satisfied. We still take $G$ to be the
movable part of $|2\sigma^*(K_{F_0})|$. Let $C$ be a generic
irreducible element of $|G|$. Similar to the situation in Case 2, it
suffices to show $|mK_{X'}||_F\supset |G|$ in order to verify
Theorem \ref{Key}(4). In fact, one has $\mathcal
{O}_{B}(2)\hookrightarrow f_*\omega_{X'}^{m_1}$ since
$n_{\Gamma}\geq 3$. Thus there is the inclusion:
$$f_*\omega_{X'/B}^2\hookrightarrow f_*\omega_{X'}^{2m_1+2}.$$
Since $f_*\omega_{X'/B}^2$ is semi-positive (= weakly positive) by
Viehweg \cite{VV} and thus generated by global sections, one has
$$m\pi^*(K_X)|_F\ge (2m_1+2)\pi^*(K_X)|_F\ge G.$$
In a word, Theorem \ref{Key} (3) and (4) are satisfied for all $m\ge
m_0+2m_1+2$.

Note that we have $p\geq 2$ by our definition. Take a $\beta$,
nearby $\frac{p}{2(m_1+p)}\geq \frac{1}{m_1+2}$, by virtue of Lemma
\ref{beta}.

If $|G|$ is composed with a pencil of curves and $g(C)=2$, then
 $\sigma^*(K_{F_0})\cdot C\ge 2$ by Lemma
\ref{10}. This gives $\xi\ge \frac{2}{m_1+2}\sigma^*(K_{F_0})\cdot
C\geq \frac{4}{m_1+2}$ by Lemma 2.8. Take $m\ge 2m_1+5$. Then
$\alpha=(m-1-\frac{m_1}{p}-\frac{1}{\beta})\xi\geq
2+\frac{4}{m_1+2}>2$. Thus Theorem \ref{Key} says that $\varphi_m$
is birational for all $m\ge \text{max}\{m_0+2m_1+2, 2m_1+5\}$.

Otherwise either $|G|$ is composed with a pencil of curves and
$g(C)\geq 3$ or $|G|$ is not composed of a pencil of curves. In the
later case, after a necessary birational modification to get the
base point freeness of $|G|$, one sees that $2g(C)-2=K_S\cdot
C+C^2\geq 4$. Again $g(C)\ge 3$. If we take a very large $m$ such
that $\alpha$ is big enough, then Theorem \ref{Key} gives
$$m\xi\ge 2g(C)-2+\roundup{(m-1-\frac{m_1}{p}-\frac{1}{\beta})\xi},$$
which means $\xi\ge \frac{8}{3m_1+6}$. Take $m\ge 3m_1+4$. Then
$\alpha\ge 2+\frac{6m_1-4}{3m_1+6}>2$. Therefore Theorem \ref{Key}
says that $\varphi_{m}$ is birational for all $m\ge
\text{max}\{m_0+2m_1+2, 3m_1+4\}$.

In a word, whenever $n_{\Gamma}\geq 3$ and $d=1$, $\varphi_{m}$ is
birational for all $m\ge \text{max}\{m_0+2m_1+2, 3m_1+4\}$.
\end{proof}

\begin{lem}\label{10} Let $S$ be a smooth projective surface of
general type with $K_{S_0}^2=1$ and $p_g(S)=0$ where $S_0$ is the
minimal model of $S$. Let $\sigma:S\longrightarrow S_0$ be the blow
down. Denote by $|C|$ the movable part of $|2K_S|$. If $g(C)=2$,
then $\sigma^*(K_{S_0})\cdot C\ge 2$.
\end{lem}
\begin{proof}
Set $\overline{C}=\sigma_*(C)$. Clearly $h^0(S_0, \overline{C})\ge
h^0(S,C)$. Thus $\overline{C}$ moves in a family. Because $|C|$ is
the movable part of $|2K_S|$, $|\overline{C}|$ must be the movable
part of $|2K_{S_0}|$ since $P_{2}(S)=P_2(S_0)$. We can write
$2K_{S_0}\sim \overline{C}+Z$ where $Z$ is the fixed part.

If $\overline{C}^2=0$, then $\overline{C}$ must be smooth and
$\sigma^*(K_{S_0})\cdot C=K_{S_0}\cdot \overline{C}\ge 2$.

If $\overline{C}^2>0$ and $K_{S_0}\cdot \overline{C}=1$, then
$\overline{C}^2\le \frac{(K_{S_0}\cdot
\overline{C})^2}{K_{S_0}^2}=1.$ Clearly $\overline{C}$ is smooth.
The Hodge index theorem says $\overline{C}\equiv K_{S_0}$. So
$Z\equiv K_{S_0}$. According to Bombieri \cite{Bom} or \cite{BPV},
$|3K_{S_0}|$ gives a birational map. So
$\varphi_{3}|_{\overline{C}}$ is birational for a general
$\overline{C}$. Because $Z\equiv K_{S_0}$ is nef and big, one has
$H^1(S_0, K_{S_0}+Z)=0$ by the vanishing theorem. So there is the
following surjective map:
$$H^0(S_0, 3K_{S_0})\longrightarrow H^0(\overline{C},
K_{\overline{C}}+Z|_{\overline{C}}).$$ Since $Z$ is effective and
$Z\cdot \overline{C}=K_{S_0}^2=1$, $Z|_{\overline{C}}$ is a single
point. So the Riemann-Roch on $\overline{C}$ gives
$h^0(K_{\overline{C}}+Z|_{\overline{C}})=2$. Thus the linear system
$|K_{\overline{C}}+Z|_{\overline{C}}|$ can only give a finite map
onto ${\mathbb P}^1$, a contradiction. Therefore $K_{S_0}\cdot
\overline{C}>1$.
\end{proof}

Theorems \ref{d3}, \ref{d2}, \ref{d1} directly imply Theorem
\ref{m0}.

Theorem \ref{main} (iv) follows from Theorem  \ref{m0}, Theorem
\ref{main} (i) and Theorem \ref{main}(ii).

\section{\bf Plurigenera of 3-folds of general type with $\chi=1$}

First let us recall Reid's plurigenus formula (at page 413 of
\cite{YPG}) for a minimal 3-fold $X$ of general type:
$$P_m(X)=\frac{1}{12}m(m-1)(2m-1)K_X^3-(2m-1)\chi(\mathcal
{O}_X)+l(m)  \eqno (3.1)$$ where $m>1$ is an integer, the correction
term
$$l(m):=\sum_{Q}l(Q,m):=\sum_{Q}\sum_{j=1}^{m-1}
\frac{\overline{bj}(r-\overline{bj})}{2r},$$ where the sum
$\sum_{Q}$ runs through all baskets Q of singularities  of type
$\frac{1}{r}(a,-a,1)$ with the positive integer $a$ coprime to $r$,
$0<a<r$, $0<b<r$, $ab\equiv 1$ (\text{mod} $r$), $\overline{bj}$ the
smallest residue of $bj$ \text{mod} $r$. Reid's result (Theorem 10.2
in \cite{YPG}) says that the above baskets $\{\text{Q}\}$ of
singularities are in fact virtual (!) and that one need not worry
about the authentic type of all those terminal singularities on $X$,
though $X$ may have non-quotient  terminal singularities.
Iano-Fletcher (\cite{Fletcher}) has showed that the set of baskets
$\{\text{Q}\}$ in Reid's formula is uniquely determined by $X$.

\begin{lem}\label{l1}\cite[Lemma 3.1]{Flt}
For all $m\ge 0$ and $n\ge 1$, one has $l(m+2n)\ge l(m)+nl(2)$ and
the equality holds if and only if all the singularities are of type
$\frac{1}{2}(1,-1,1)$.
\end{lem}

\begin{lem}\label{l2}\cite[lemma 3.2]{Flt}
For positive integers $\alpha>\beta>n$,
$$l(\frac{1}{\alpha}(1,-1,1),n)\ge l(\frac{1}{\beta}(1,-1,1),n).$$
\end{lem}

\begin{lem}\label{l3}\cite[lemma 3.3]{Flt} For an arbitrary positive integer
 $a$ coprime to $r$ and for a positive integer
 $n\le [\frac{r+1}{2}]$,
$$l(\frac{1}{r}(a,-a,1),n)\ge l(\frac{1}{r}(1,-1,1),n)$$ where
$[\frac{r+1}{2}]$ denotes the integral part of
$\frac{r+1}{2}$.
\end{lem}

\begin{lem}\label{l4}\cite[corollary 3.4]{Flt}
For positive integers $a$, $\alpha$, $\beta$ with $0\le \beta \le
\alpha$ and $a$ coprime to $\alpha$ and for a positive integer $n\le
[\frac{\alpha+1}{2}]$,
$$l(\frac{1}{\alpha}(a,-a,1),n)\ge l(\frac{1}{\beta}(1,-1,1),n).$$
\end{lem}

\begin{setup}\label{assumption}{\bf Assumption}. {}From now on within
this section we assume $X$ to be a minimal projective 3-fold of
general type with only ${\mathbb Q}$-factorial terminal
singularities and with $\chi(\mathcal {O}_X)=1$.
\end{setup}

\begin{prop}\label{p2}
 If $P_{2}(X)\ge 1$, then $P_{2n}(X)\ge n$ and $P_{2n+1}(X)\ge n-1$ for
all integers $n\ge 2$.
\end{prop}
\begin{proof}
Since
\begin{eqnarray*}
&P_{4}(X)=7K_{X}^{3}-7+l(4),\\
&P_{2}(X)=\frac{1}{2}K_{X}^{3}-3+l(2).
\end{eqnarray*}
Lemma \ref{l1} says $l(4)\ge 2 l(2)$. Thus
$$P_{4}(X)-2P_{2}(X)=6K_{X}^{3}-1+(l(4)-2l(2))$$
which implies $P_{4}(X)\ge 6K_{X}^{3}+1$. So one has $P_{4}(X)\ge
2.$

Assume $P_{2n}(X)\ge n$ for any integer $n\ge 2$. One has
\begin{eqnarray*}
&&P_{2n+2}(X)=\frac{(2n+2)(2n+1)(4n+3)}{12}K_{X}^{3}-(4n+3)+l(2n+2),\\
&&P_{2n}(X)=\frac{2n(2n-1)(4n-1)}{12}K_{X}^{3}-(4n-1)+l(2n),\\
&&P_{2n+2}(X)-P_{2n}(X)-P_{2}(X)=kK_{X}^{3}-1+l(2n+2)-l(2n)-l(2)
\end{eqnarray*}\
where $k>0$ and $l(2n+2)-l(2n)-l(2)\ge 0$ by Lemma \ref{l1}. Thus
$P_{2n+2}(X)\ge n+kK_{X}^{3}>n$ which implies $P_{2n+2}(X)\ge n+1$.
The first assertion is proved.

Now we study $P_{2n+1}(X)$. Similarly one has
\begin{eqnarray*}
&P_{3}(X)=\frac{5}{2}K_{X}^{3}-5+l(3),\\
&P_{5}(X)=15K_{X}^{3}-9+l(5),\\
&P_{5}(X)-P_{3}(X)-P_{2}(X)=12K_{X}^{3}-1+l(5)-l(3)-l(2).
\end{eqnarray*}
So $P_{5}(X)\ge P_{3}(X)+12K_{X}^{3}$ which says $P_{5}(X)\ge 1$.
Assume that $P_{2n+1}(X)\ge n-1$ for a number $n\ge 2$. Then a
calculation gives:
\begin{eqnarray*}
&&P_{2n+3}(X)-P_{2n+1}(X)-P_{2}(X)\\
&=&tK_{X}^{3}-1+l(2n+3)-l(2n+1)-l(2)
\end{eqnarray*}where $t>0$. Thus $P_{2n+3}(X)\ge n$.
We are done.
\end{proof}


\begin{prop}\label{s28} Assume $P_{2}(X)=0$.
If $X$ contains a virtue basket Q of singularities with canonical
index $r(\text{Q})\geq 28$, Then $P_{18}(X)\ge 3$.
\end{prop}
\begin{proof} If $X$ contains a virtue basket $Q$ with index
$r(\text{Q})=r\ge 37$, Lemma \ref{l1}, Lemma \ref{l2} and Lemma
\ref{l3} give
\begin{eqnarray*}
&&l(18)\ge l(Q,18)\ge l(\frac{1}{r}(1,-1,1),18)\\
&= &\sum_{j=1}^{17}\frac{j(r-j)}{2r}\ge
\sum_{j=1}^{17}\frac{j(37-j)}{74}>52.
\end{eqnarray*} Thus
$P_{18}(X)>-35+l(18)>17$.

If $X$ contains a virtue basket $Q$ with index $r\in [28,36] $, one
can verify $l(Q,18)>37$ case by case. So $P_{18}(X)>2$.
Alternatively, one may apply the property of the polynomial
$y=x(r-x)$ to greatly simply the calculation.
\end{proof}

\begin{thm}\label{18} Let $X$ be a minimal projective 3-fold of
general type with ${\mathbb Q}$-factorial terminal singularities
and with $\chi(\mathcal {O}_X)=1$. Then $P_{2l+8}(X)\ge 3$ for all
$l\ge 0$.
\end{thm}
\begin{proof}
If $P_2(X)\ge 1$, then Proposition \ref{p2} implies $P_{18}(X)\ge
9$. Assume $P_2(X)=0$ from now on. Proposition \ref{s28} tells that
we may even assume the index $r(Q)\le 27$ for all virtue basket Q of
$X$. This makes it possible for us to study within limited
possibilities. The table in the last part lists all possible types
of Q with index $\le 27$.\medskip

{\bf Step 1}. $P_{18}(X)\ge 1$.

To the contrary, assume $P_{18}=0$. Then
$P_{3}(X)=P_{6}(X)=P_{9}(X)=0$. For a positive integer $n$, set
$$\Delta_n:=n^2l(2)+l(n)-l(n+1).$$
One has
\begin{eqnarray*}
&&P_{3}(X)-5P_{2}(X)=10-\sum_{Q}\Delta_{2},\\
&&P_{6}(X)-P_{3}(X)-50P_{2}(X)=144-\sum_{Q}(\Delta_{3}+
    \Delta_{4}+\Delta_{5}),\\
&&P_{9}(X)-P_{6}(X)-149P_{2}(X)=441-\sum_{Q}(\Delta_{6}+
    \Delta_{7}+\Delta_{8}),\\
&&P_{18}(X)-P_{9}(X)-1581P_{2}(X)=4725-\sum_{Q}(\Delta_{9}+\cdots
+\Delta_{17}).
\end{eqnarray*}

Set $F$ to be the matrix:
$$ \left(
\begin{array}{cccc}
1 &-11 &0 &-4\\
0 &1 &-3 &-11\\
0 &0 &1 &-7\\
0 &0 &0 &1\\
\end{array}
\right)
$$

Since $P_{3}(X)=P_{6}(X)=P_{9}(X)=P_2(X)=0$, we get
$$ \left\{
\begin{array}{l}
\nabla_1:=\sum_{Q}\Delta_{2}=10\\
\nabla_2:=\sum_{Q}(\Delta_{3}+\Delta_{4}+\Delta_{5})=144\\
\nabla_3:=\sum_{Q}(\Delta_{6}+\Delta_{7}+\Delta_{8})=441\\
\nabla_4:=\sum_{Q}(\Delta_{9}+\cdots+\Delta_{17})=4725\\
\end{array}
\right.
$$
Then
\begin{eqnarray*}
(\nabla_1',\nabla_2',\nabla_3',\nabla_4')&=&
(\nabla_1,\nabla_2,\nabla_3,\nabla_4)F\\
&=&(10,34,9,14). \end{eqnarray*}

On the other hand, for any basket Q, we can formally compute
$\Delta_i(Q)$ for any positive integer $i$. So one gets
$\nabla_j(Q)$ for $j=1,\cdots,4$. Taking the product with the matrix
$F$, one even gets $\nabla_j'(Q)$. In our table, we have listed all
those values of $\nabla_j'(Q)$ where $j=1,\cdots,4$. Searching with
a computer, one finds that the only possible combinations of baskets
Q of $X$ are as the following:
\begin{enumerate}

\item[(i)] 3 of type $\frac{1}{2}(1,-1,1)$, 2 of type $\frac{1}{5}(3,-3,1)$ and
1 of type $\frac{1}{10}(7,-7,1)$;

\item[(ii)] 4 of type $\frac{1}{2}(1,-1,1)$, 3 of type $\frac{1}{3}(2,-2,1)$, 1
of type $\frac{1}{5}(4,-4,1)$ and 1 of type $\frac{1}{5}(3,-3,1)$;

\item[(iii)] 2 of type $\frac{1}{2}(1,-1,1)$, 2 of type $\frac{1}{3}(2,-2,1)$,
1 of type $\frac{1}{4}(3,-3,1)$ and 1 of type
$\frac{1}{12}(7,-7,1)$.
\end{enumerate}
For each case, one gets $l(2)=3$ which means $K_X^3=0$ a
contradiction to $X$ being of general type. Therefore $P_{18}(X)\ge
1$.\medskip

{\bf Step 2}. $P_{18}(X)\ge 2$.

Similarly we assume $P_{18}=1$. Then there are 5 possibilities:

\begin{enumerate}
\item[(a)] $P_{3}(X)=P_{6}(X)=P_{9}(X)=0$;

\item[(b)] $P_{3}(X)=P_{9}(X)=0$, and $P_{6}(X)=1$;

\item[(c)] $P_{3}(X)=P_{6}(X)=0$, and $P_{9}(X)=1$;

\item[(d)] $P_{3}(X)=P_{6}(X)=P_{9}(X)=1$;

\item[(e)] $P_{3}(X)=0$ and $P_{6}(X)=P_{9}(X)=1$.
\end{enumerate}

 In the case (a), one has
\begin{eqnarray*}
(\nabla_1',\nabla_2',\nabla_3',\nabla_4')&=&
(\nabla_1,\nabla_2,\nabla_3,\nabla_4)F\\
&=&(10,34,9,13). \end{eqnarray*} Searching with a computer, one
finds that the only possible combination of baskets Q of $X$ is:

(iv) 5 of type $\frac{1}{2}(1,1,1)$, 2 of type $\frac{1}{4}(3,1,1)$,
1 of type $\frac{1}{5}(4,1,1)$ and 1 of type $\frac{1}{5}(3,2,1)$.

A calculation shows $l(2)=3$. Then $K_{X}^{3}=0$ which contradicts
to $X$ being of general type.

In the case (b) through (e), one has, respectively:
\begin{eqnarray*}
&&(\nabla_1',\nabla_2',\nabla_3',\nabla_4')=
(\nabla_1,\nabla_2,\nabla_3,\nabla_4)F\\
&=&(10,33,13,17), (10, 34,8,21), (9,45,9,18),(10,33,12,25).
\end{eqnarray*}
Searching with a computer, one finds that the only possible
combinations of baskets Q of $X$ occurring in case (b) and case (d)
are:
\begin{enumerate}
\item[(v)] (case b) 5 of type $\frac{1}{2}(1,1,1)$, 4 of type $\frac{1}{3}(2,1,1)$
and 1 of type $\frac{1}{6}(5,1,1)$; (One gets $l(2)=3$.)

\item[(vi)] (case d)  2 of type $\frac{1}{4}(3,1,1)$, 1 of type $\frac{1}{5}(3,2,1)$,
1 of type $\frac{1}{7}(4,3,1)$ and 1 of type $\frac{1}{8}(5,3,1)$;
(One gets $l(2)>3$.)

\item[(vii)] (case d) 9 of type $\frac{1}{3}(2,1,1)$. (One gets $l(2)=3$.)

\end{enumerate}
Clearly one obtains $K_{X}^{3}\le 0$, a contradiction to $X$ being
of general type. Therefore $P_{18}(X)\ge 2$. \medskip

{\bf Step 3}. $P_{18}(X)\ge 3$.

Assume $P_{18}(X)=2$. There are still five possibilities (a) through
(e) as listed above for $P_3(X)$, $P_6(X)$ and $P_9(X)$. Then one
gets corresponding datum as follows:
\begin{eqnarray*}
(\nabla'_1,\nabla'_2,\nabla'_3,\nabla'_4)& =&(10,34,9,12),
(10,33,13,16), (10,34,8,20)\\
&& (9,45,9,17), (10,33,12,24).
\end{eqnarray*}
Still searching with a computer, one gets possible combinations of
baskets Q of $X$:

\begin{enumerate}

\item[(viii)] (Case a) 2 of type $\frac{1}{2}(1,1,1)$, 2 of type $\frac{1}{3}(2,1,1)$, 1
of type $\frac{1}{4}(3,1,1)$, 1 of type $\frac{1}{5}(3,2,1)$ and 1
of type $\frac{1}{7}(5,2,1)$;

\item[(ix)] (Case b) 6 of type $\frac{1}{2}(1,1,1)$, 1 of type $\frac{1}{3}(2,1,1)$, 2
of type $\frac{1}{4}(3,1,1)$ and 1 of type $\frac{1}{6}(5,1,1)$;

\item[(x)]  (Case c) 4 of type $\frac{1}{2}(1,1,1)$, 2 of type $\frac{1}{5}(4,1,1)$
and 2 of $\frac{1}{5}(3,2,1)$;

\item[(xi)]  (Case d) 1 of type $\frac{1}{2}(1,1,1)$, 6 of type $\frac{1}{3}(2,1,1)$
and 2 of type $\frac{1}{4}(3,1,1)$.
\end{enumerate}
For the situation (viii), one obtains
$P_{2}(X)=P_{3}(X)=P_{4}(X)=\cdots=P_{11}(X)=0$, $P_{12}(X)=1$,
$P_{13}(X)=0$, $P_{14}(X)=P_{15}(X)=P_{16}(X)=P_{17}(X)=1$,
$P_{18}(X)=P_{19}(X)=2$, $P_{20}(X)=P_{21}(X)=3$.
\medskip

\noindent{\bf Claim}. {\em Situation (viii) doesn't exist.}
\begin{proof}
According to Reid (see (10.3) of \cite{YPG}), one has
$$\frac{1}{12}K_X\cdot c_2(X)=-2\chi({\mathcal
O}_X)+\sum_Q\frac{r_Q^2-1}{12r_Q}$$ where $c_2(X)$ is defined via
the intersection theory by taking a resolution of singularities
over $X$. Miyaoka (Corollary 6.7 of \cite{Miyaoka}) says $K_X\cdot
c_2(X)\geq 0$. Thus one sees the inequality
$$\sum_Q\frac{r_Q^2-1}{r_Q}\geq 24\chi({\mathcal
O}_X).$$ Now since the datum of (viii) doesn't fit into the above
inequality, situation (viii) doesn't exist at all.
\end{proof}

For other situations (ix) through (xi), one gets $l(2)=3$. Thus
$K_X^3=0$ which is impossible.

In a word, we have proved $P_{18}(X)\ge 3$.
\medskip

{\bf Step 4}. $P_{18+2l}(X)\ge 3$ for all $l\geq 0$.

Since
$$P_{20}(X)=\frac{14820}{12}K_{X}^{3}-39+l(20),$$
$$P_{20}(X)-P_{18}(X)-P_{2}(X)=qK_{X}^{3}-1+l(20)-l(18)-l(2)\ge
-1+qK_{X}^{3}$$ where $q>0$. So $P_{20}(X)\ge P_{18}(X)\ge 3$.

Assume $P_{18+2k}\ge 3$ for any integer $k\ge 1$. One has
$$P_{2k+20}(X)=q_1K_{X}^{3}+(-4k-39)+l(2k+20),$$
$$P_{2k+18}(X)=q_2K_{X}^{3}+(-4k-35)+l(2k+18),$$
$$P_{2k+20}(X)-P_{2k+18}(X)-P_{2}(X)=q_3K_{X}^{3}-1+l(2k+20)-l(2k+18)-l(2)$$ where
$q_1$, $q_2$, $q_3>0$ and $l(2k+20)-l(2k+18)-l(2)\geq 0$ by Lemma
\ref{l1}. Thus $P_{2k+20}(X)\ge 2+qK_{X}^{3}>2$ which implies
$P_{2k+20}(X)\ge 3$. We are done.
\end{proof}

\begin{thm}\label{14}  Let $X$ be a minimal projective 3-fold of
general type with ${\mathbb Q}$-factorial terminal singularities and
with $\chi(\mathcal {O}_X)=1$. Then $P_{m}(X)\ge 1$ for all $m\ge
14$.
\end{thm}
\begin{proof}  First we show $P_{2n}(X)\ge 1$ for all $n\ge 7$.
In fact, $P_{14}(X)-P_{12}(X)-P_{2}(X)>-1+l(14)-l(12)-l(2)$ and
$P_{12}(X)\ge 1$ by Fletcher \cite{Flt}. Thus  $P_{14}(X)\ge 1$.
Assume that $P_{2n}(X)\ge 1$ for some $n\ge 7$. Then
$P_{2n+2}(X)-P_{2n}(X)-P_{2}(X)>-1+l(2n+2)-l(2n)-l(2)$ which
implies $P_{2n+2}(X)> P_{2n}(X)-1\ge 0$.

Next we assume $P_{15}(X)\ge 1$. With a similar method, one can
see $P_{2n+1}(X)\ge 1$ for all $n\ge 8$. Actually
$P_{2n+1}(X)-P_{2n-1}(X)-P_{2}(X)>-1+l(2n+1)-l(2n-1)-l(2)$ which
implies $P_{2n+1}(X)> P_{2n-1}(X)-1\ge 0$. So $P_{2n+1}(X)\ge 1$
whenever $P_{2n-1}(X)\ge 1$.

Now we consider what happens when $P_{15}(X)=0$. Clearly
$P_{3}(X)=P_{5}(X)=0$. By Proposition \ref{p2}, we may assume
$P_2(X)=0$.

Clearly we have
$$ \left\{
\begin{array}{l}
P_{3}(X)-5P_{2}(X)=10-\sum_{Q}\Delta_{2}(Q),\\
P_{5}(X)-P_{3}(X)-25P_{2}(X)=71-\sum_{Q}(\Delta_{3}(Q)+\Delta_{4}(Q)),\\
P_{15}(X)-P_{5}(X)-985P_{2}(X)=2935-\sum_{Q}(\Delta_{5}(Q)+\cdots+
\Delta_{14}(Q)).
\end{array}
\right.
$$
Then
$$ \left\{
\begin{array}{l}
\sum_{Q}\Delta_{2}(Q)=10\\
\sum_{Q}(\Delta_{3}(Q)+\Delta_{4}(Q))=71\\
\sum_{Q}(\Delta_{5}(Q)+\cdots+\Delta_{14}(Q))=2935\\
\end{array}
\right.
$$

\begin{claim}
If X contains a virtue basket Q of type $\frac{1}{r}(a,-a,1)$ where
$r\ge 26$, then $P_{15}(X)\ge 1$.
\end{claim}
\begin{proof}
According to Lemmas \ref{l2}, \ref{l3}, \ref{l4}, one has
$$P_{15}(X)> -29+l(\frac{1}{29}(1,-1,1), 15)\ge 6$$
whenever $r\ge 29$. For $26\le r\le 28$, one may check, case by
case, the inequality $P_{15}(X)\ge 1$.
\end{proof}

Therefore it suffices to consider all those virtue basket Q with
$r\le 25$. Set $\Lambda_{1}=\Delta_{2}$,
$\Lambda_{2}=\Delta_{3}+\Delta_{4}$ and
$\Lambda_{3}=\Delta_{5}+\cdots+\Delta_{14}$. Set $G$ to be the
following matrix:
$$ \left(
\begin{array}{rrr}
1 &-5 &-5\\
0 &1 &-40\\
0 &0 &1\\
\end{array}
\right)
$$
Let
$$(\Lambda_1', \Lambda_2',\Lambda_3') =(\Lambda_1,
\Lambda_2,\Lambda_3)G =(10,21,45). $$ Again we can search with a
computer. The only possibilities of combinations of
 virtue baskets Q of $X$ are:

\begin{enumerate}

\item[(xii)] 1 of type $\frac{1}{10}(7,3,1)$, 2 of type
$\frac{1}{5}(3,2,1)$ and 3 of type $\frac{1}{2}(1,1,1)$; ($l(2)=3$)

\item[(xiii)] 1 of type $\frac{1}{4}(3,1,1)$, 2 of type
$\frac{1}{8}(5,3,1)$ and 3 of type $\frac{1}{2}(1,1,1)$; ($l(2)=3$)

\item[(xiv)] 2 of type $\frac{1}{3}(2,1,1)$, 1 of type
$\frac{1}{4}(3,1,1)$, 1 of type $\frac{1}{12}(7,5,1)$ and 2 of type
$\frac{1}{2}(1,1,1)$; ($l(2)>3$)

\item[(xv)] 3 of type $\frac{1}{3}(2,1,1)$, 1 of type
$\frac{1}{5}(4,1,1)$, 1 of type $\frac{1}{5}(3,2,1)$ and 4 of type
$\frac{1}{2}(1,1,1)$; ($l(2)=3$)

\item[(xvi)] 4 of type $\frac{1}{3}(2,1,1)$, 1 of type
$\frac{1}{6}(5,1,1)$ and 5 of type $\frac{1}{2}(1,1,1)$. ($l(2)=3$)
\end{enumerate}

Clearly $K_{X}^{3}\le 0$, a contradiction. Therefore, there are no
3-folds of general type with $\chi(\mathcal {O})=1$ and
$P_{15}(X)=0$. We have showed $P_{15}(X)>0$ which implies the whole
theorem.
\end{proof}

\begin{proof}[{\bf Proof of Theorem \ref{main}(iii)}]
By Theorem \ref{18}, we have $P_{18}\geq 3$. So we may take
$m_1=18$. Consider the map $\varphi_{18}$ and keep the same notation
as in \ref{setup} and \ref{Key}.

If $d=3$, then $\varphi_m$ is birational for all $m\geq 56$ by
Theorem \ref{d3}.

If $d=2$, since $m_1=18$, $\varphi_m$ is birational for all $m\geq
63$ by Theorem \ref{d2}.

If $d=1$, then $\varphi_m$ is birational for all $m\geq 58$ by
Theorem \ref{d1}(3).

So Theorem \ref{main}(iii) is proved.
\end{proof}

\newpage

\renewcommand{\arraystretch}{1.5}

$$
\begin{array}{ccccccccc}
No. \hspace*{2mm}&\hspace*{2mm}Singularity\hspace*{2mm} &\hspace*{2mm}\nabla'_{1} \hspace*{2mm}&\hspace*{2mm}\nabla'_{2} \hspace*{2mm}&\nabla'_{3} \hspace*{2mm}&\hspace*{2mm}\nabla'_{4}\hspace*{2mm} &\hspace*{2mm}\Lambda'_{1} \hspace*{2mm}&\hspace*{2mm}\Lambda'_{2} \hspace*{2mm}&\Lambda'_{3}\\
1 &\frac{1}{2}(1,-1,1)  &1 &1 &1 &-1       &1 & 1 &0\\
2 &\frac{1}{3}(2,-2,1)  &1 &5 &1 &2        &1 &3 & 1\\
3 &\frac{1}{4}(3,-3,1)  &1 &7 &1 &3        &1  & 4 &1\\
4 &\frac{1}{5}(4,-4,1)  &1 &8 &1 &10       &1  &  4 & 25\\
5 &\frac{1}{5}(3,-3,1)  &2 &7 &1 &2        &2  &  4  & 17\\
6 &\frac{1}{6}(5,-5,1)  &1 &8 &4 &14       &1  &  4  & 41\\
7 &\frac{1}{7}(6,-6,1)  &1 &8 &6 &18       &1  &  4  & 52\\
8 &\frac{1}{7}(5,-5,1)  &3 &8 &3 &5        &3  &  5 & 24\\
9 &\frac{1}{7}(4,-4,1)  &2 &12 &3 &7       &2  &  7  &  8\\
10 &\frac{1}{8}(7,-7,1) &1 &8 &7 &25       &1  &  4  & 59\\
11 &\frac{1}{8}(5,-5,1) &3 &12 &3 &3       &3  &  7 & 22\\
12 &\frac{1}{9}(8,-8,1) &1 &8 &7 &35       &1 & 4 & 65\\
13 &\frac{1}{9}(7,-7,1) &4 &9 &4 &10       &4 & 6 &27\\
14 &\frac{1}{9}(5,-5,1) &2 &15 &2 &17      &2 & 8 &29\\
15 &\frac{1}{10}(9,-9,1)  &1 &8 &7 &43     &1 & 4 &70\\
16 &\frac{1}{10}(7,-7,1)  &3 &17 &4 &13    &3  & 10 & 11\\
17 &\frac{1}{11}(10,-10,1) &1 &8 &7 &50     &1 & 4 & 74\\
18 &\frac{1}{11}(9,-9,1)  &5 &10 &5& 13    &5 & 7 &29\\
19 &\frac{1}{11}(8,-8,1)  &4 &17 &4 &8     &4 & 10 &25\\
20 &\frac{1}{11}(7,-7,1)  &3 &19 &4 &12    &3 & 11 &10\\
21 &\frac{1}{11}(6,-6,1)  &2 &16 &5 &27    &2 & 8 &67\\
22 &\frac{1}{12}(11,-11,1) &1 &8 &7 &56     &1 & 4&77\\
23 &\frac{1}{12}(7,-7,1)  &5 &15 &4 &9     &5 & 9 &42\\
24 &\frac{1}{13}(12,-12,1) &1 &8 &7 &61     &1  &  4  &79\\
25 &\frac{1}{13}(11,-11,1) &6 &11 &6 &15    &6 & 8 &30\\
26 &\frac{1}{13}(10,-10,1) &4 &22 &5 &17    &4  & 13& 13\\
27 &\frac{1}{13}(9,-9,1)  &3 &22 &3 &22    &3 & 12 &31\\
28 &\frac{1}{13}(8,-8,1)  &5 &19 &4 &6     &5 & 11 &40\\
29 &\frac{1}{13}(7,-7,1)  &2 &16 &10 &33   &2 & 8 &94
\end{array}
$$
$$
\begin{array}{ccccccccc}
No. \hspace*{2mm}&\hspace*{2mm}Singularity\hspace*{2mm} &\hspace*{2mm}\nabla'_{1} \hspace*{2mm}&\hspace*{2mm}\nabla'_{2} \hspace*{2mm}&\nabla'_{3} \hspace*{2mm}&\hspace*{2mm}\nabla'_{4}\hspace*{2mm} &\hspace*{2mm}\Lambda'_{1} \hspace*{2mm}&\hspace*{2mm}\Lambda'_{2} \hspace*{2mm}&\Lambda'_{3}\\
30 &\frac{1}{14}(13,-13,1)  &1 &8 &7 &65          &1 & 4 &80\\
31 &\frac{1}{14}(11,-11,1)  &5 &22 &5 &12         &5 & 13 &27\\
32 &\frac{1}{14}(9,-9,1)   &3  &23 &3 &28        &3  & 12  &55\\
33 &\frac{1}{15}(14,-14,1)  &1 &8 &7 &68          &1 & 4 &80\\
34 &\frac{1}{15}(13,-13,1)  &7 &12 &7 &16         &7 & 9 &30\\
35 &\frac{1}{15}(11,-11,1)  &4 &26 &5 &16         &4 & 15 &11\\
36 &\frac{1}{15}(8,-8,1)   &2 &16 &13 &   44     &2 & 8&111\\
37 &\frac{1}{16}(15,-15,1)& 1 & 8   &  7   & 70   &1 & 4 & 80\\
38 &\frac{1}{16}(13,-13,1)  & 5&27&6 & 20         &5 & 16 & 14\\
39 &\frac{1}{16}(11,-11,1)  & 3  & 24 & 6 & 38    &3 & 12 &92\\
40 &\frac{1}{16}(9,-9,1)   & 7 & 17 & 7 & 16     &7  & 11 &51\\
41 &\frac{1}{17}(16,-16,1)  & 1 & 8 & 7 & 71      &1 & 4 &80\\
42 &\frac{1}{17}(15,-15,1)  & 8 & 13 & 8 & 16     &8 & 10 &30\\
43 &\frac{1}{17}(14,-14,1)  & 6 & 27 & 6 & 15     &6 & 16 &28\\
44 &\frac{1}{17}(13,-13,1)  & 4 & 29 & 4 & 26     &4 & 16 &32\\
45 &\frac{1}{17}(12,-12,1)  & 7 & 22 & 5 & 12     &7 & 13 &59\\
46 &\frac{1}{17}(11,-11,1)  & 3 & 24 & 9 & 42     &3 & 12&108\\
47 &\frac{1}{17}(10,-10,1)  & 5 & 29 & 7 & 21     &5  & 17 & 19\\
48 &\frac{1}{17}(9,-9,1)   & 2 & 16 & 14 & 61    &2 & 8 &124\\
49 &\frac{1}{18}(17,-17,1)  & 1 & 8 & 7 & 71      &1 & 4 &80\\
50 &\frac{1}{18}(13,-13,1)  & 7 & 26 & 5 & 8      &7 & 15 &57\\
51  &\frac{1}{18}(11,-11,1) & 5 & 31 & 7 & 19     &5 & 18 &18\\
52  &\frac{1}{19}(18,-18,1) & 1 & 8 & 7 & 71      &1 & 4 &80\\
53  &\frac{1}{19}(17,-17,1) & 9 & 14 & 9 & 15     &9 & 11 &30\\
54  &\frac{1}{19}(16,-16,1) & 6 & 32 & 7 & 22     &6  & 19&15\\
55  &\frac{1}{19}(15,-15,1) & 5 & 33 & 6 & 19     &5 & 19 &12\\
56  &\frac{1}{19}(14,-14,1) & 4 & 31 & 4 & 38     &4 & 16 & 80\\
57  &\frac{1}{19}(13,-13,1) & 3 & 24 & 14 & 47    &3 & 12&135\\
58  &\frac{1}{19}(12,-12,1) & 8 & 23 & 7 & 14     &8 & 14 &66\\
59  &\frac{1}{19}(11,-11,1) & 7 & 29 & 7 & 11     &7 & 17 &47\\
\end{array}
$$
$$
\begin{array}{ccccccccc}
No. \hspace*{2mm}&\hspace*{2mm}Singularity\hspace*{2mm} &\hspace*{2mm}\nabla'_{1} \hspace*{2mm}&\hspace*{2mm}\nabla'_{2} \hspace*{2mm}&\nabla'_{3} \hspace*{2mm}&\hspace*{2mm}\nabla'_{4}\hspace*{2mm} &\hspace*{2mm}\Lambda'_{1} \hspace*{2mm}&\hspace*{2mm}\Lambda'_{2} \hspace*{2mm}&\Lambda'_{3}\\
60  &\frac{1}{19}(10,-10,1)  & 2 & 16 & 14 & 78            &2 & 8&135\\
61  &\frac{1}{20}(19,-19,1)  & 1 & 8 & 7 & 71              &1  &  4 &80\\
62  &\frac{1}{20}(17,-17,1)  & 7 & 32 & 7 & 17             &7 & 19 &29\\
63  &\frac{1}{20}(13,-13,1)  & 3 & 24 & 16 & 51            &3 & 12&146\\
64  &\frac{1}{20}(11,-11,1)  & 9 & 19 & 9 & 23             &9 & 13 &56\\
65 &\frac{1}{21}(20,-10,1)   & 1 & 8 & 7 & 71              &1 & 4 &80\\
66  &\frac{1}{21}(19,-19,1)  & 10 & 15 & 10 & 14           &10 & 12 &30\\
67 &\frac{1}{21}(17,-17,1)   & 5 & 36 & 5 & 29             &5 & 20 &33\\
68 &\frac{1}{21}(16,-16,1)   & 4 & 32 & 7 & 48             &4 & 16&117\\
69  &\frac{1}{21}(13,-13,1)  & 8 & 31 & 7 & 9              &8  & 18 &62\\
70 &\frac{1}{21}(11,-11,1)  & 2 & 16 & 14 & 93            &2 & 8&144\\
71 &\frac{1}{22}(21,-21,1)   & 1& 8 & 7 &71                &1 & 4 &80\\
72 &\frac{1}{22}(19,-19,1)   & 7 & 37 & 8 &24              &7 & 22 &16\\
73 &\frac{1}{22}(17,-17,1)   & 9 &29 &6&14                 &9 & 17 &76\\
74&\frac{1}{22}(15,-15,1) & 3 & 24&   19  &62              &3 & 12&163\\
75&\frac{1}{22}(13,-13,1) &5& 37 &5 &39                    &5 & 20 &60\\
76&\frac{1}{23}(22,-22,1) & 1 & 8 & 7 &71                  &1 & 4 &80\\
77&\frac{1}{23}(21,-21,1) & 11 &16 &11 &13                 &11 & 13 &30\\
78&\frac{1}{23}(20,-20,1)    & 8&  37    &    8    &19     &8 & 22 &30\\
79&\frac{1}{23}(19,-19,1)    & 6&  40  &  7 &22            &6 & 23 &13\\
80&\frac{1}{23}(18,-18,1)    & 9&  33    &    6    &10     &9 & 19 &74\\
81&\frac{1}{23}(17,-17,1)    & 4&  32    &   13  &56       &4 & 16&149\\
82&\frac{1}{23}(16,-16,1)    &10  &  25  & 10 &21          &10 & 16 &75\\
83&\frac{1}{23}(15,-15,1)     & 3&  24  & 20  &69          &3 & 12&170\\
84&\frac{1}{23}(14,-14,1)     & 5&  38 & 5 &45             &5 & 20 &84\\
85&\frac{1}{23}(13,-13,1)    & 7&  39    &    9    &30    &7 & 23 &24\\
86&\frac{1}{23}(12,-12,1)    & 2&  16    &   14    &106   &2 & 8&151\\
87&\frac{1}{24}(23,-23,1)     & 1&   8    &    7    &71    &1 & 4 &80\\
88&\frac{1}{24}(19,-19,1)     & 5&  39    &    5    &48     &5  & 20&105\\
89&\frac{1}{24}(17,-17,1)     & 7&  41    &   10    &28     &7 & 24 &27\\
\end{array}
$$
$$
\begin{array}{ccccccccc}
No. \hspace*{2mm}&\hspace*{2mm}Singularity\hspace*{2mm} &\hspace*{2mm}\nabla'_{1} \hspace*{2mm}&\hspace*{2mm}\nabla'_{2} \hspace*{2mm}&\nabla'_{3} \hspace*{2mm}&\hspace*{2mm}\nabla'_{4}\hspace*{2mm} &\hspace*{2mm}\Lambda'_{1} \hspace*{2mm}&\hspace*{2mm}\Lambda'_{2} \hspace*{2mm}&\Lambda'_{3}\\
90&\frac{1}{24}(13,-13,1)   & 11 &21 &   11 &28           &11 & 15 &59\\
91&\frac{1}{25}(24,-24,1)  & 1& 8 & 7  &71                &1 & 4 &80\\
92&\frac{1}{25}(23,-23,1)   &12     &  17    &   12 &12   &12 & 14 &30\\
93&\frac{1}{25}(22,-22,1)    & 8&  42    &    9    &26    &8 & 25 &17\\
94&\frac{1}{25}(21,-21,1)   & 6&  43    &    6    &32     &6 & 24 &34\\
95&\frac{1}{25}(19,-19,1)  & 4&  32    &   18    &61      &4 & 16&176\\
96&\frac{1}{25}(18,-18,1)  & 7&  43    &   10    &26      &7 & 25 &26\\
97&\frac{1}{25}(17,-17,1)   & 3&  24    &   21    &86     &3 & 12&183\\
98&\frac{1}{25}(16,-16,1)    &11     &  26    &   11 &26  &11 &  17 &78\\
99&\frac{1}{25}(14,-14,1)  & 9&  39    &    9    &20     &9 & 23 &52\\
100&\frac{1}{25}(13,-13,1)  & 2&  16    &   14    &117    &2  &  8  &156\\
101&\frac{1}{26}(25,-25,1)  & 1&   8    &    7    &71    \\
102&\frac{1}{26}(23,-23,1) & 9&  42    &    9    &21      \\
103&\frac{1}{26}(21,-21,1)  & 5 &40    &    8    &58        \\
104&\frac{1}{26}(19,-19,1)& 11 & 31 & 10 &19                  \\
105&\frac{1}{26}(17,-17,1) & 3 & 24 & 21 &96                    \\
106&\frac{1}{26}(15,-15,1) & 7 & 45 & 9 &28                      \\
107&\frac{1}{27}(26,-26,1) & 1 & 8 & 7 &71\\
108&\frac{1}{27}(25,-25,1) & 13 & 18 & 13 &11\\
109&\frac{1}{27}(23,-23,1) & 7 & 47 & 8 &25\\
110&\frac{1}{27}(22,-22,1) & 11 & 36 & 7 &16 \\
111&\frac{1}{27}(20,-20,1)& 4 & 32 & 22 &69 \\
112&\frac{1}{27}(19,-19,1) & 10 & 41 & 10 &14 \\
113&\frac{1}{27}(17,-17,1) & 8 & 46 & 11 &34\\
114&\frac{1}{27}(16,-16,1) & 5 & 40 & 11 &65\\
115&\frac{1}{27}(14,-14,1)& 2 & 16 & 14 &126\\
\end{array}
$$

\begin{rem} Almost one year after the first version of this paper was put
to arXiv, Yongnam Lee informed us of the relevant paper \cite{shin}.
We admit that the effectivity of Miyaoka-Reid inequality (see the
proof of Claim) was first observed in \cite{shin}. In fact, one may
take $m_0=6$ and $m_1=10$ by virtue of \cite{shin}. Thus Theorem
\ref{m0} imply the following:
\end{rem}

\begin{cor} Let $V$ be a nonsingular projective 3-fold of general
type with $\chi({\mathcal O}_V)=1$. Then $\varphi_m$ is birational
onto its image for all $m\geq 54$.
\end{cor}


\end{document}